\documentclass{sigcomm-alternate}
\usepackage{amsfonts,latexsym,amsmath,amstext,amssymb,verbatim,epsfig,psfrag}
\usepackage{colordvi,pstcol}
\usepackage{graphicx,psboxit}
\usepackage{psfrag,graphicx,balance}

\newtheorem{remark}{Remark}

\begin{document}
\title{\ttlit D-iteration: Evaluation of the Asynchronous Distributed Computation}

\numberofauthors{1}
\author{
   \alignauthor Dohy Hong\vspace{2mm}\\
   \affaddr{Alcatel-Lucent Bell Labs}\\
   \affaddr{Route de Villejust}\\
   \affaddr{91620 Nozay, France}\\
   \email{\normalsize dohy.hong@alcatel-lucent.com}
}

\date{\today}
\maketitle

\begin{abstract}
The aim of this paper is to present a first evaluation of the potential of an asynchronous distributed computation associated to the recently proposed approach, D-iteration: the D-iteration is a fluid diffusion based iterative method, which has the advantage of being natively distributive. It exploits a simple intuitive decomposition of the matrix-vector product as elementary operations of fluid diffusion associated to a new algebraic representation.
We show through experiments on real datasets how much this approach can improve the computation efficiency when the parallelism is applied: with the proposed solution, when the computation is distributed over $K$ virtual machines (PIDs), the memory size to be handled by each virtual machine decreases linearly with $K$ and the computation speed increases almost linearly with $K$ with a slope becoming closer to one when the number $N$ of linear equations to be solved increases.
 
\end{abstract}
\category{G.1.0}{Mathematics of Computing}{Numerical Analysis}[Parallel algorithms]
\category{G.1.3}{Mathematics of Computing}{Numerical Analysis}[Numerical Linear Algebra]
\category{C.2.4}{Computer Systems Organization}{Computer-Communication Networks}[Distributed Systems]
\terms{Algorithms, Performance}
\keywords{Distributed computation, Iteration, Fixed point, Eigenvector.}
\begin{psfrags}
\section{Introduction}\label{sec:intro}
Surprisingly enough, there was a recent result \cite{dohy} showing a potentially significant
reduction of the numerical computation cost in a very classical problem of
the calculation of the eigenvector of a large sparse matrix, based on a new representation/interpretation/decomposition of
the matrix-vector product as elementary operations of fluid diffusion (cf. \cite{d-algo, serge}).
This is an alternative solution to existing iterative methods
(cf. \cite{Golub1996, Saad, Bagnara95aunified}): its potential in the context of
PageRank equation has been shown in \cite{dohy} and the application of the approach
in a general context is described in \cite{d-algo} (D-iteration).

The complexity of the computation of the eigenvector of a matrix is a very well
known problem and it increases rapidly with the dimension of the vector space.
Efficient, accurate methods to compute eigenvectors of arbitrary matrices
are in general a difficult problem (cf. power iteration \cite{mises}, QR algorithm \cite{francis, kub}).

Starting from the algorithm proposed in \cite{dohy}
(a computation speed-up by factor 4-40 depending on the graphs was observed), we
present here a first evaluation of the speed-up factor that can be cumulatively
applied on the previous performance when a distributive architecture is used.
There have been a lot of researches concerning the distributed computation of the
linear equations (cf. \cite{Bertsekas:1989:PDC:59912, jela, Lubachevsky:1986:CAA:4904.4801, DBLP:journals/corr/abs-cs-0606047, Kohlschutter06efficientparallel}),
with a particular interest on asynchronous iteration scheme.
The elimination of the synchronization constraints is important on heterogeneous platforms
for a better efficiency and for an effective scalability to distributed platforms.
The distributive algorithm defined and evaluated here is directly inspired from the first asynchronous distributed scheme 
that has been described in \cite{distributed}.

We recall that the PageRank equation can be written under the form:
\begin{eqnarray}\label{eq:affine}
X = P' X = P.X + B
\end{eqnarray}
where $P$ and $P'$ are non-negative square matrices of size $N\times N$
($P'$ is the stochastic matrix associated to the transition probability), 
$B$ a non-negative vector of size $N$ and $X$ the eigenvector
of $P'$ ($x_i$ gives the score of the importance of the page $i$).

We recall that the D-iteration approach works when the spectral radius
of $P$ is strictly less than 1 (so that the power series $\sum_{n\ge 0} P^n B$ is convergent).
In this paper, we will focus only on the
specific case when $P$ is strictly contractive: for the PageRank equation, 
the damping factor $d$ is the explicit contraction coefficient (the usually
considered value of $d$ is $0.85$).

The main objective of this paper is to evaluate through simple simulations (based on a single PC
for now) the efficiency of the D-iteration based simple distributed computation architecture,
when applied on a large matrix (current simulations limited to $N=10^6$) in the context of PageRank type equation.

In Section \ref{sec:algo}, we define the distributed algorithm that is used. 
The evaluation of the computation efficiency in the context of the PageRank
equation is studied in Section \ref{sec:eval}.

\section{Distributed algorithm}\label{sec:algo}
The fluid diffusion model in the general case is described
by the matrix $P$ associated with a weighted graph ($p_{ij}$ is
the weight of the edge from $j$ to $i$, positive or negative) and the initial condition
$F_0 = B$ (cf. \cite{d-algo}).

We recall the definition of the two vectors used in the D-iteration method:
the fluid vector $F_n$ is defined by:
\begin{eqnarray}
F_n &=& (I_d - J_{i_n} + P J_{i_n}) F_{n-1}.\label{eq:defF}
\end{eqnarray}
where:
\begin{itemize}
\item $I_d$ is the identity matrix;
\item $I = \{i_1, i_2, ..., i_n,...\}$ with $i_n \in \{1,..,N\}$ is a deterministic or random sequence
  such that the number of occurrence of each value $k\in \{1,..,N\}$ in $I$ is infinity;
\item $J_k$ a matrix with all entries equal to zero except for
  the $k$-th diagonal term: $(J_k)_{kk} = 1$.
\end{itemize}

And the history vector $H_n$ defined by ($H_0$ initialized to a null vector):
\begin{eqnarray}\label{eq:defH}
H_n &=& \sum_{k=1}^n J_{i_k} F_{k-1}.
\end{eqnarray}

The distance to the limit $X$ ($L_1$ norm $|X-H_n|$) is then bounded by
$r/(1-d)$, where $r$ is the residual fluid $r = |F_n|$ and
where $d$ is the contraction factor of $P$
(equal to the damping factor for PageRank equation).
Note that $H_n$ is an increasing function for each component (when $P$ and $B$ are non-negative), so that:
$|X-H_n| = \sum_{i=1}^N |x_i - (H_n)_i| = \sum_{i=1}^N (x_i - (H_n)_i)$
Therefore, $|X-H_n| = 1 - |H_n|$ when $X$ is a probability vector.

\subsection{Pseudo-code}
We recall that the above equations \eqref{eq:defF} and \eqref{eq:defH} are the
mathematical formulations of the following algorithm (D-iteration):

\begin{verbatim}
Initialization:
  H[i] := 0;
  F[i] := B_i;
  r := |F|;

Iteration:
k := 1;
While ( r/(1-d) > Target_Error )
  Choose i_k;
  sent := F[i_k];
  H[i_k] += sent;
  F[i_k] := 0;
  If ( i_k has at least one child )
    For all child node j of i_k:
      F[j] += sent * p(j,i_k);
  r := |F|;
  k++;
\end{verbatim}
\begin{remark}
Of course, we don't have to compute the quantity $r$ in each step of iteration. 
\end{remark}
\subsection{PageRank equation}
We recall that in the context of PageRank equation, $P$ is of the form:
$dQ$ where $Q$ is a stochastic matrix (or sub-stochastic matrix in presence
of dangling nodes) and satisfies the equation:
\begin{eqnarray}\label{eq:pr}
 X = d Q X + (1-d) V,
\end{eqnarray}
with $V$ a personalization vector (by default, equal to the uniform probability vector $(1/N,..,1/N)^t$)
and $d$ is the damping factor.
Equivalently, we can write for each entry of $X$:
\begin{eqnarray}\label{eq:pr2}
 x_i = \frac{(1-d)}{N} + d \sum_j q_{ij} x_j.
\end{eqnarray}

And we have, $q_{ij} = \frac{1}{\#out_j}$, where
$\#out_j$ is the number of outgoing links from $j$ (when there is at least one
outgoing link from $j$, otherwise, one may complete the $j$-th column of $Q$ by $1/N$
cf. \cite{deep, Boldi2009, dohy}). 
\subsection{Partition sets}
In the following, we consider two simple $K$ partitioning sets:
\begin{itemize}
\item Uniform partition: $\Omega_1 =\{1, 2, ..., N/K\}$,
  $\Omega_2 =\{N/K+1, 2, ..., 2\times N/K\}$, etc
\item Cost Balanced (CB) partition: $\Omega_k = \{\omega_k, \omega_k+1, ..., \omega_{k+1}-1\}$ such
  that $\sum_{n=\omega_k}^{\omega_{k+1}-1} (\#out_n) = L/K$, where $L$ is the total number of non-null
  entries of $Q$ (or the total number of links of the graph associated to $Q$),
\end{itemize}
such that $\{1,...,N\} = \Omega = \cup_k \Omega_k$.
The intuition of the cost balanced partition is that when we apply the diffusion iteration on
all nodes of each $\Omega_k$, the diffusion cost is constant.
The main reason why we chose this is the simplicity of its computation: 
as a possibly better alternative solution, we could equalize the
quantity: $\sum_{n,m = \omega_k}^{\omega_{k+1}-1} p_{n,m}$ (which gives a very rough information on
at which speed the residual fluid - cf. Section \ref{sec:pid} - of the set $\Omega_k$ disappears).

As described in \cite{distributed}, each $PID_k$ computes $(F_n)_i$ and $(H_n)_i$ for $i\in\Omega_k$
($[F_n]_k$ and $[H_n]_k$)
and those informations are exchanged between PIDs.
In the following, we assumed that there is no communication cost when information is exchanged between PIDs.

\begin{remark}
We believe that there is here some nice simple adaptation scheme to be designed making the set $\Omega_k$
dynamically evolve in time, with the idea that an idle PID could take nodes from the busiest
PID: such an operation is indeed quite easy with our diffusion model: we just need to transmit
the information $F_n$ and $H_n$ of the candidate nodes (then, possibly re-forwarding fluid to the
new destination PID, if not well synchronized etc).
\end{remark}
\subsection{PID modelling}\label{sec:pid}
We consider a time stepped approximation for the simulation of the distributed computation cost
(for now running on a single PC):
during each time step, each PID can execute $PID\_Speed_k$ operations.
By default, we set: $PID\_Speed_k = PID\_Speed = L/K$ (PIDs are assumed to compute at the same speed).

We define an elementary operation cost as a diffusion cost from one node to another node.

In the following, the number of iterations is defined as the normalized elementary operation cost dividing
the number of operations by $L$ (so that it can be easily compared to the cost of one matrix-vector product,
or one iteration in power iteration).

Each $PID_k$ applies the diffusion from nodes in $\Omega_k$: the diffusion iteration is done
until the residual fluid 
$$
r_k = |[F_n]_k| = \sum_{i\in\Omega_k}(F_n)_i
$$ 
of $\Omega_k$ is below a certain threshold:
$$
r_k < T.
$$
When this target is reached, this $PID_k$ can go in sleep mode.

In each $\Omega_k$, the choice of the node sequence $I_k$ is done based on the ideas
proposed in \cite{dohy, d-algo}:
$$i_n^{(k)} = \arg\max_{i\in\Omega_k} \left((F_{n-1})_i/((\#in_i +1) \times (\#out_i+1))\right)$$
with $\#in_i$ (resp. $\#out_i$) equal to the number of incoming (resp. outgoing) 
links to (resp. from) the node $i$.
To approximate the computation of the $\arg\max$ (the computation of the exact value is
computation costly and not necessary),
we used a cyclic test on nodes and all $i$ such that $(F_n)_i$ is above
a threshold $T_k$ are chosen, then we scale down the threshold progressively
by a factor $\alpha>1$:
\begin{verbatim}
  T_k \= alpha.
\end{verbatim}
We observed that the results below are not sensitive to the choice of $\alpha$
(between 1.2 and 2, by default we chose 1.5).
\subsection{Fluid transmission to external nodes}
By external nodes, we mean here all nodes $j\notin\Omega_k$ such that $PID_k$
has $(F_n)_j$ waiting for the transmission to a $PID_{k'}$ ($j\in\Omega_{k'}$).
Each PID monitors its residual fluid $r_k$ and the fluid amount to be sent $s_k$.
The quantity $s_k$ can be computed in two ways:
\begin{itemize}
\item continuously from regular updates when diffusions are applied in $\Omega_k$;
\item periodically, by the computation of the increment of $H$ (from the last time the fluids
  have been sent to external nodes) and applying diffusion on this increment (cf. \cite{distributed}).
\end{itemize}

This last procedure is preferred, because in this way, we can save the computation cost, by applying
(note that the quantity $r_k$ need not to be updated in every computation step as well) the following
algorithm (*):
\begin{verbatim}
  Choose i (in Omega_k);
  sent := F[i];
  F[i] := 0;
  H[i] += sent;
  If ( i has at least one child in Omega_k )
    For all child node j of i (in Omega_k):
      F[j] += sent * p(j,i);
  k++;
\end{verbatim}

At the time fluids are sent to external nodes we store $[H_n]_k$ in \verb+H_old+, then for the next external
diffusion, we can apply $\left(P(H-H_{old})\right)_{i\notin \Omega_k}$ to find the fluid amount
to be sent to each external node $i$.

The transmission of fluid $F_n$ from $PID_k$ is done when the total fluid to be transmitted
is above certain threshold: by default, we used the condition:
\begin{eqnarray}\label{eq:send}
s_k &>& r_k/K.
\end{eqnarray}

In fact, we can even further optimize the information exchange on $(F_n)_{j\notin\Omega_k}$: this
can be done directly by the receivers and the sender $PID_k$ can only send $[H_n]_k = (H_n)_{i\in\Omega_k}$.
\begin{itemize}
\item each $PID_k$ computes $[F_n]_k$ and $[H_n]_k$ based on (*);
\item when the condition (such as) \eqref{eq:send} is satisfied, $PID_k$ sends $[H_n]_k$ to all
  other PIDs (we could also centralize the reception of this to a dedicated PID).
\end{itemize}

The drawback of this approach is that the receiver $PID_{k'}$ needs to compute the quantity 
$\left(P([H_n]_k-[H_{old}]_k)\right)_{i\in \Omega_{k'}}$ ($[H_{old}]_k$ is the lastly received $[H]_k$ from $PID_k$),
which means that a vector of size $N$ needs to be stored and instead of storing only the column vectors of $P$
corresponding to the set $\Omega_k$:
$(P)_{i\in\Omega, j\in\Omega_k}$, we have to store 
$(P)_{i\in\Omega_k, j\in\Omega}$
(cf. \cite{distributed}).

\subsection{Fluid reception from external nodes}
When new fluids are received from other PIDs, we need two operations:
\begin{itemize}
\item add them to existing fluids (node per node; can be delayed);
\item update the threshold $T_k$, rescaling up by a factor: in our experiments, we applied a factor proportional to
  $\min((r_k + received)/r_k, received)$ (if $r_k > 0$, otherwise, $T_k$ is set to $received$).
\end{itemize}

\section{Evaluation of the distributed computation cost}\label{sec:eval}
The aim of this section is to show the potential of the distributed computation gain
when the D-iteration is used.

For the evaluation purpose, we experimented the D-iteration on a
web graph imported from the dataset \verb+uk-2007-05@1000000+
(available on \verb+http://law.dsi.unimi.it/datasets.php+) which has
41.247.159 links on 1.000.000 nodes (45.766 dangling nodes).

Below we vary $N$ from 1000 to 1000000, extracting from the dataset the
information on the first $N$ nodes.

\begin{table}
\begin{center}
\begin{tabular}{|l|ccc|}
\hline
N & L (nb links) & L/N & D (Nb dangling nodes)\\
\hline
1000 & 12.935 & 12.9 & 41 (4.1\%)\\
10000 & 125.439 & 12.5 & 80 (0.8\%)\\
100000 & 3.141.476 & 31.4 & 2729 (2.7\%)\\
1000000 & 41.247.159 & 41.2 & 45766 (4.6\%)\\
\hline
\end{tabular}\caption{Extracted graph: $N=1000$ to $1000000$.}
\end{center}
\end{table}


\subsection{Uniform partition}
We first analyse the results obtained with a uniform partition.
Figure \ref{fig:2PIDs-1000} shows a typical result with $K=2$: we compare the evolution
of the cost in terms of the number of normalized iteration for a targeted 
error (an upper bound of the distance to the limit: the y-axis shows the value $r_k/(1-d)$ for each PID). 
In this case, PID1 
computes the first half nodes and the last half nodes is handled by PID2.
Because of the possible idle states of PIDs, at the end of each test, the computation
cost of different PIDs are not necessarily equal. And naturally the global speed of the
convergence is dictated by the slowest one. We observe that the average
cost of PID1 and PID2 for a given error (horizontal line) is roughly half
of the cost with a single PID. 
When $N=100000$, we see that the work load is
not well balanced: the slopes are different almost
by factor 2.

\begin{figure}[htbp]
\centering
\includegraphics[angle=-90, width=\linewidth]{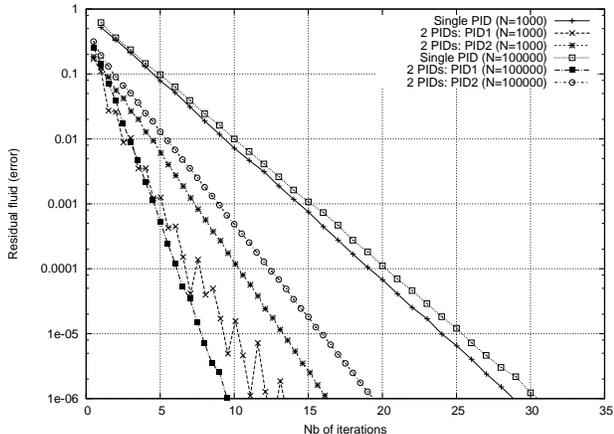}
\caption{$N=1000$ and $N=100000$: 1 PID vs 2 PIDs.}
\label{fig:2PIDs-1000}
\end{figure}

Figure \ref{fig:4PIDs-1000} shows the same result with 4 PIDs. Of course, when
partitioned uniformly there is no reason that the slope of each PID is the same:
in this case, the convergence speeds of the different PIDs are very different.
\begin{figure}[htbp]
\centering
\includegraphics[angle=-90, width=\linewidth]{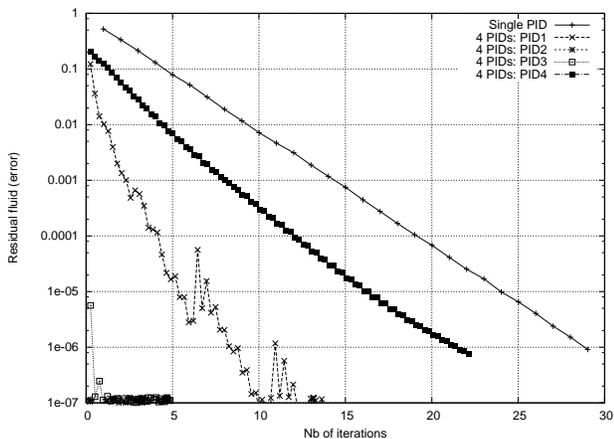}
\caption{$N=1000$: 1 PID vs 4 PIDs.}
\label{fig:4PIDs-1000}
\end{figure}

Note that because of the computation cost reduction applied for the external diffusion,
each PID's slope is in fact closely linked to the spectral radius of the extracted square matrix
corresponding to the nodes of the set $\Omega_k$ (similar to $d$ which plays a role for a single PID case).

Hence, to optimize the distributed computation efficiency, there are two aspects:
\begin{itemize}
\item first, the slope of the convergence should be equalized between PIDs;
\item second, there should be no idle PIDs.
\end{itemize}
If we don't consider the dynamic evolution of the partition sets,
the above second point can not be controlled.
However, the first point can be approximately optimized with the CB partition (cf. Section \ref{sec:cb}).
Our current understanding and guess is that this CB partition is not too bad
for a very large matrix. This will be illustrated below.

Figure \ref{fig:4PIDs-1000-fluid} is a visualization of the residual fluid $r_k$
and the cumulated fluid for transmission $s_k$: such a view is important to check
the stability of the algorithm (for instance, we should not have $s_k$ too far above
$r_k$).
 
\begin{figure}[htbp]
\centering
\includegraphics[angle=-90, width=\linewidth]{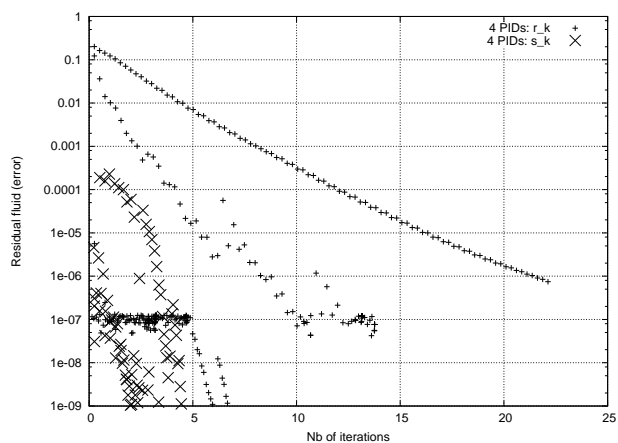}
\caption{$N=1000$: 4 PIDs: $r_k$ and $s_k$.}
\label{fig:4PIDs-1000-fluid}
\end{figure}

Now, we evaluate the distributed computation gain for different values of $N$ and $K$:
here the convergence time is arbitrarily defined as the normalized number of iterations 
(given by the slowest PID) necessary to reach a target error of $1/N$.
The results are summarized in Figure \ref{fig:cost-UP}: when $K$ is too large, there is
less gain (for small $N$). 
This was expected, because if the sizes of $\Omega_k$ are too small,
there is no possible improvement. Also when $N/K$ is too small, the information exchange cost between
PIDs becomes higher (which is not taken into account here).
However, we see clearly that when $N$ is large, the value
of $K$ on which the gain stays linear is larger. And this is the property that makes our
approach promising. From this figure, we can possibly expect to build a real implementation
on hundreds of virtual machines  that could deliver a gain over a factor of 100 
for say $N > 10^6$ (in the figure, it is reached for $N=100000$, $K=1024$).

\begin{figure}[htbp]
\centering
\includegraphics[angle=-90, width=\linewidth]{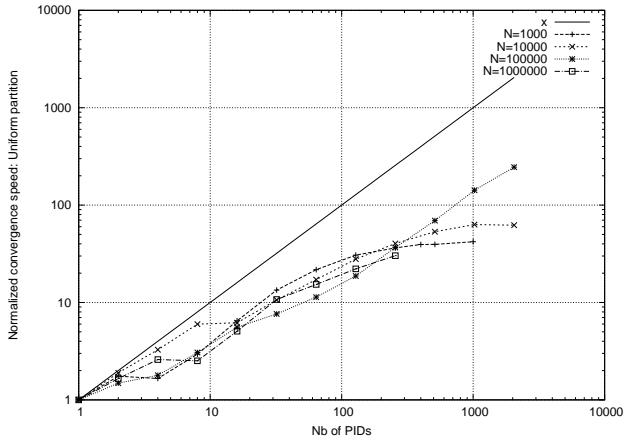}
\caption{Computation speed-up factor with $K$: uniform partition.}
\label{fig:cost-UP}
\end{figure}

\subsection{CB partition}\label{sec:cb}
The illustration of the impact of CB partition is shown in Figure \ref{fig:2PIDs-100000-CB}
($|\Omega_1|=68671$, $|\Omega_2|=31329$):
this is the typical impact that was generally observed: we see that the slopes of PIDs with
CB partition are much closer than before. This is what we expected, but this is of course not
always the case (when $N/K$ is not large enough).
Figure \ref{fig:4PIDs-100000-compa} (the respective values of $\omega_k$ are 0, 42110, 68671, 86399, 100000) 
shows another case where the slopes are made closer.

\begin{figure}[htbp]
\centering
\includegraphics[angle=-90, width=\linewidth]{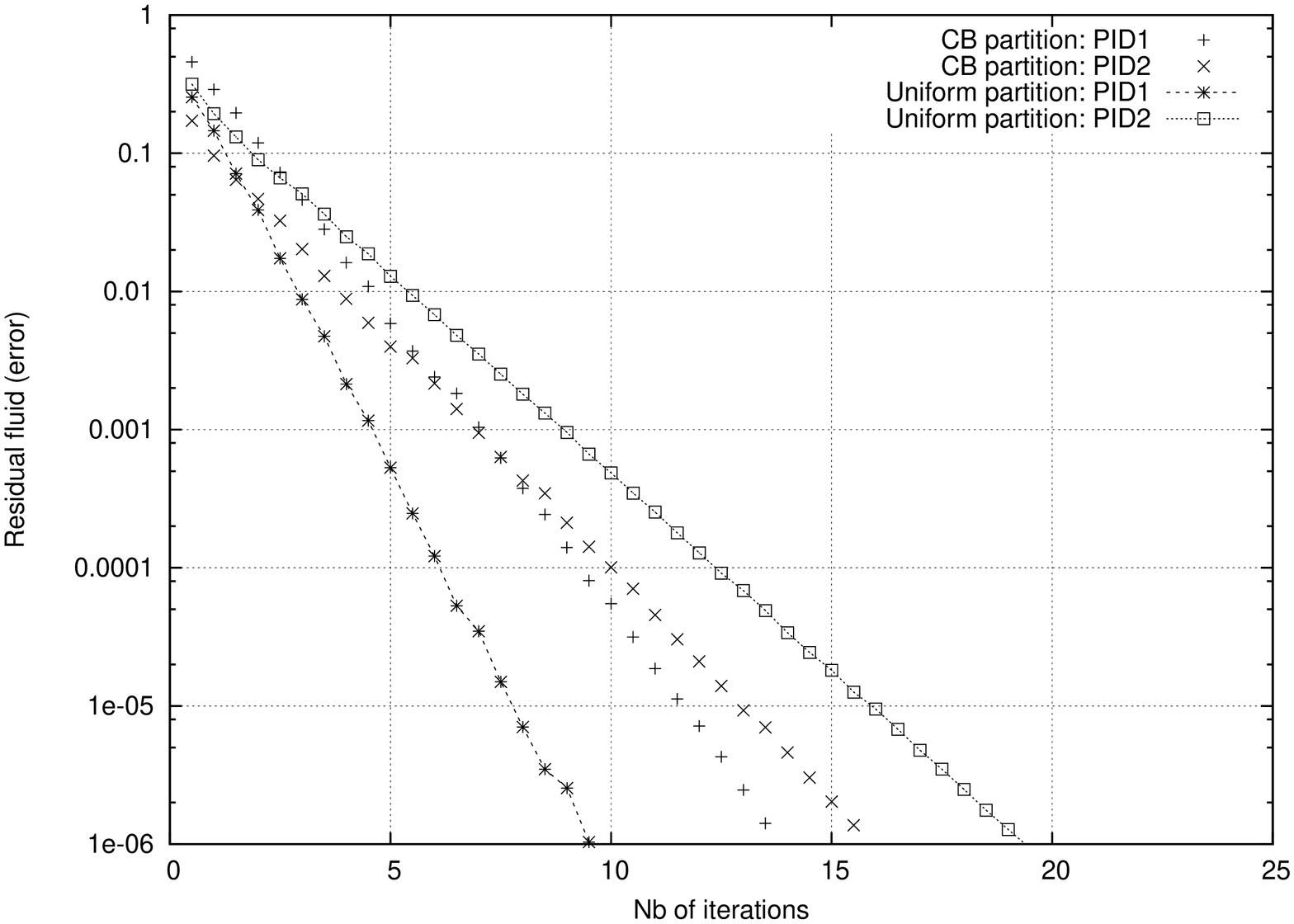}
\caption{$N=100000$: impact of CB partition.}
\label{fig:2PIDs-100000-CB}
\end{figure}

\begin{figure}[htbp]
\centering
\includegraphics[angle=-90, width=\linewidth]{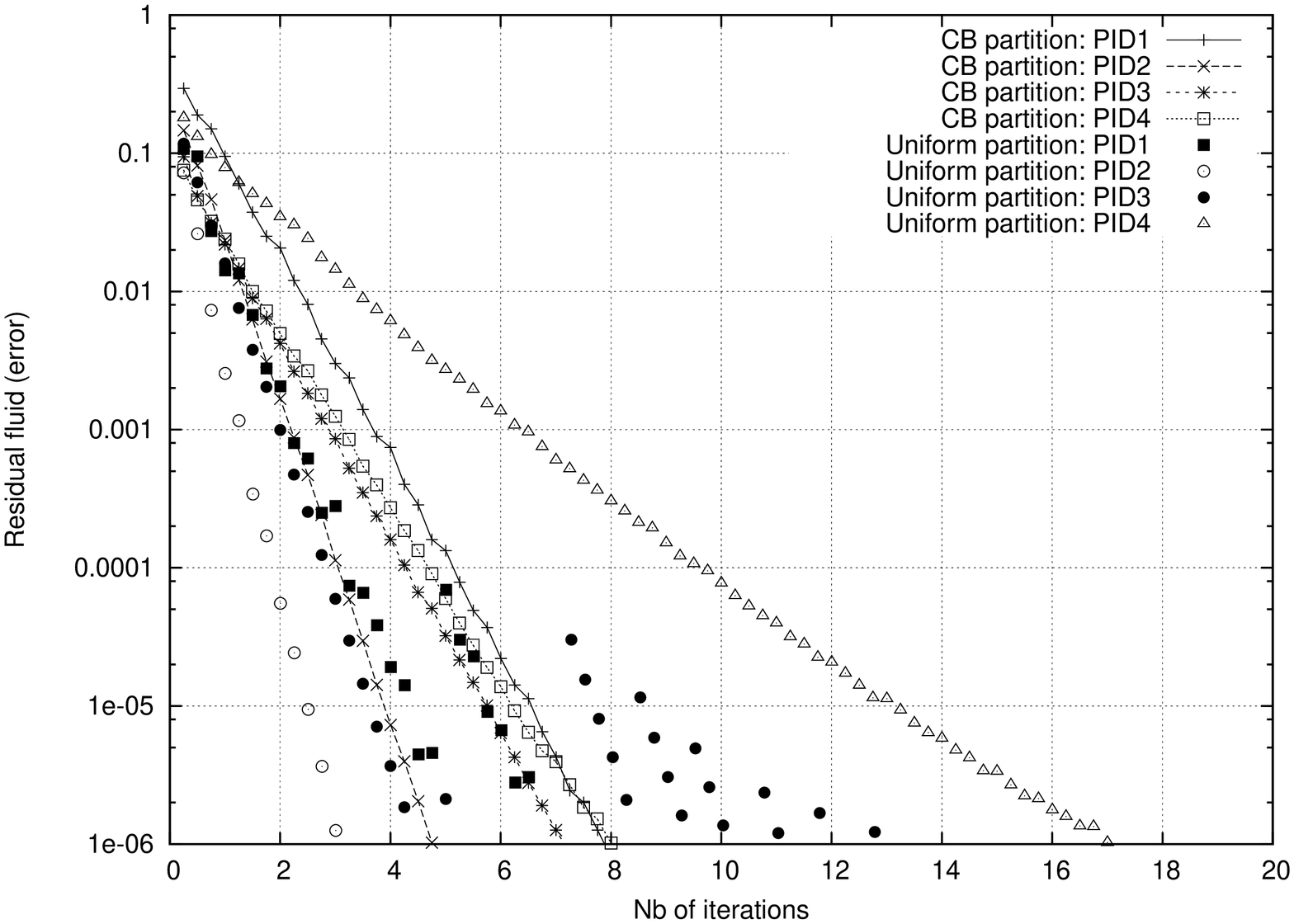}
\caption{$N=100000$: impact of CB partition.}
\label{fig:4PIDs-100000-compa}
\end{figure}

As before, the distributed computation gain for different values of $N$ and $K$
are summarized in Figure \ref{fig:cost-CB}: we see that generally the uniform partition case
is improved: Figure \ref{fig:cost-compa} shows the gain in percentage when CB partition is
used instead of uniform partition: the gain depends very much on the initial graph
configuration (we may have naturally {\em well balanced} matrix or the opposite): 
we observed here up to 250\% gain. When $K$ is too large, CB becomes meaningless.

\begin{figure}[htbp]
\centering
\includegraphics[angle=-90, width=\linewidth]{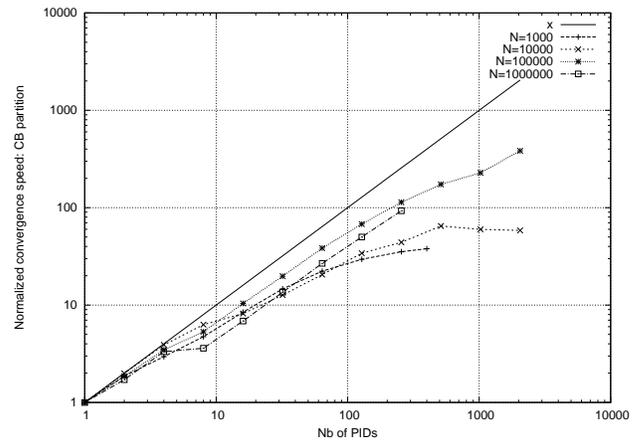}
\caption{Computation speed-up factor with $K$: CB partition.}
\label{fig:cost-CB}
\end{figure}

\begin{figure}[htbp]
\centering
\includegraphics[angle=-90, width=\linewidth]{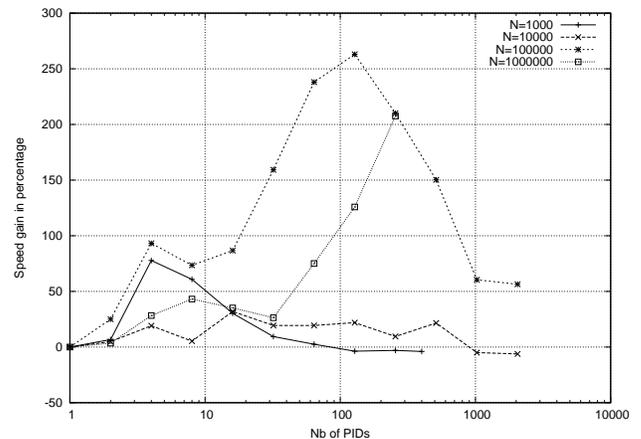}
\caption{Speed-up gain CB/Unif.}
\label{fig:cost-compa}
\end{figure}

\begin{remark}
In Figure \ref{fig:cost-compa}, the missing points are due to the impossibility
(can not get $K>N$) or the limitation of the memory size on a single PC (for $N=1000000$, $K$ up to 256).
\end{remark}

\subsection{Information exchange cost}
In this section, we want to evaluate the impact of the delayed information exchange
(for instance, due to network delays): for that purpose, in the simulation scheme above
we added a parameter \verb+Delay_Proba+ such that when each PID decides to exchange
information the transmission is tested every $PID\_Speed$ steps and delayed with probability
\verb+Delay_Proba+: the transmission is with 0 delay with probability $1-Delay\_Proba$ 
and the average transmission delay is 
$$
PID\_Speed\times Delay\_Proba/(1-Delay\_Proba).
$$

We recall that the convergence time was defined for the target error value of
$1/N$. For the required number of iterations, the normalized operation cost
includes here the elementary operation cost and the idle time 
(when $PID\_Speed$ is not totally consumed before
entering in the idle state, the remaining {\em capacity} is counted as idle time).

\begin{figure}[htbp]
\centering
\includegraphics[angle=-90, width=\linewidth]{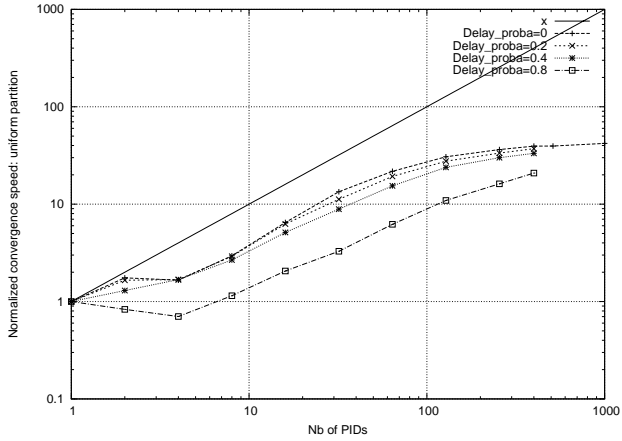}
\caption{Uniform partition, $N=1000$: impact of delayed information exchange.}
\label{fig:cost-delay-UP-1000}
\end{figure}

\begin{figure}[htbp]
\centering
\includegraphics[angle=-90, width=\linewidth]{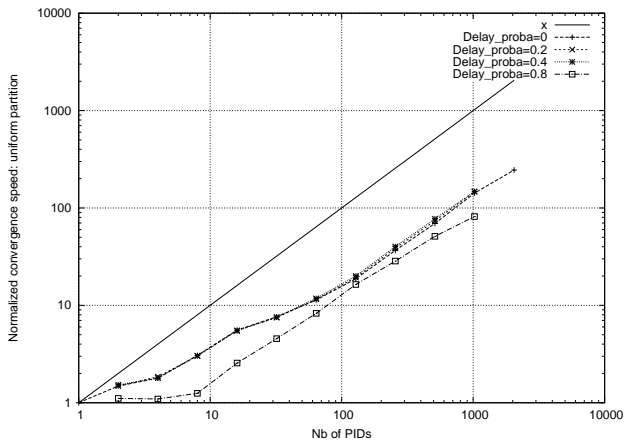}
\caption{Uniform partition, $N=100000$: impact of delayed information exchange.}
\label{fig:cost-delay-UP-100000}
\end{figure}

The results are shown in Figure \ref{fig:cost-delay-UP-1000} and Figure \ref{fig:cost-delay-UP-100000}
for $N=1000$ and $N=100000$ when uniform partition is applied: we observe that the impact is
very limited except when $Delay\_Proba$ is equal to 0.8.

Figure \ref{fig:cost-delay-CB-1000} and Figure \ref{fig:cost-delay-CB-100000}
for $N=1000$ and $N=100000$ show the impact when CB partition is applied: we observe again a very limited
impact when $Delay\_Proba$ is not too large.

\begin{figure}[htbp]
\centering
\includegraphics[angle=-90, width=\linewidth]{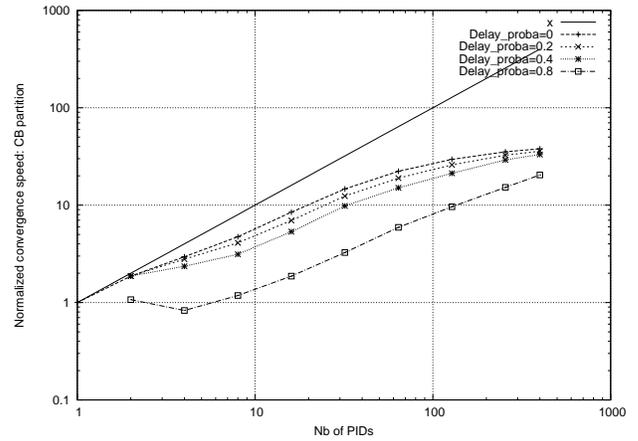}
\caption{CB partition, $N=1000$: impact of delayed information exchange.}
\label{fig:cost-delay-CB-1000}
\end{figure}

\begin{figure}[htbp]
\centering
\includegraphics[angle=-90, width=\linewidth]{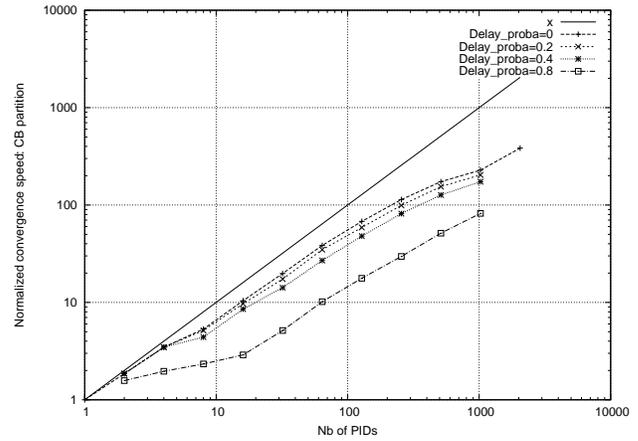}
\caption{CB partition, $N=100000$: impact of delayed information exchange.}
\label{fig:cost-delay-CB-100000}
\end{figure}

Figure \ref{fig:cost-delay-1000} and Figure \ref{fig:cost-delay-100000}
show the proportion of the idle states observed in these simulations.
The proportion of the idle state is naturally defined as:
$$
idle\_times/(elementary\_cost + idle\_times).
$$
As expected, we observe that when $K$ increases, the proportion of
the idle state is increased (this is very clear for $N=100000$ and $K>10$).

\begin{figure}[htbp]
\centering
\includegraphics[angle=-90, width=\linewidth]{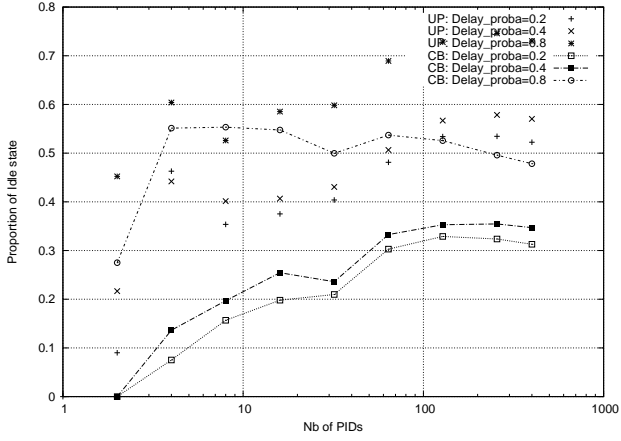}
\caption{$N=1000$: proportion of the idle state.}
\label{fig:cost-delay-1000}
\end{figure}

\begin{figure}[htbp]
\centering
\includegraphics[angle=-90, width=\linewidth]{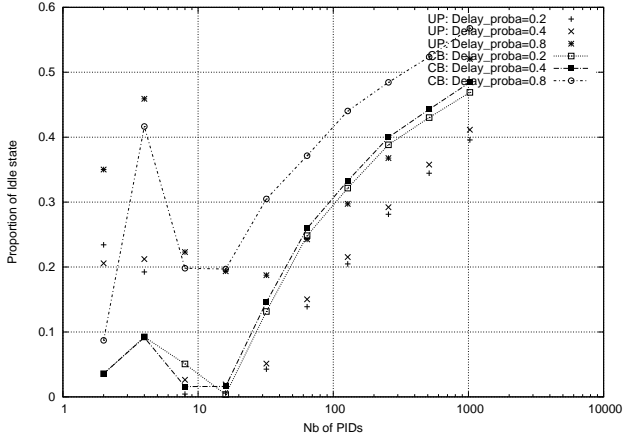}
\caption{$N=100000$: proportion of the idle state.}
\label{fig:cost-delay-100000}
\end{figure}

The above results clearly shows that the impact of the delayed information is very limited
when the main computation cost is in the PIDs internal calculation cost (and the tolerance
is quite large).
This means that thanks to our distributed computation scheme that is really asynchronous,
we can obtain a substantial gain with $K$ PIDs as soon as the diffusion computation time
on $\Omega_k$ is not significantly below the information exchange time. 

\subsection{A very simple test of the partition set adaptation}
We chose here a case when the uniform partition did not correctly balance PIDs' work load:
$N=100000$ and $K=2$ .
Figure \ref{fig:adap} gives an illustration of a very simple dynamic adaptation scheme:
in this case, we dynamically adjusted the value $\omega_2$ ($+/- 10$\%), 
starting with $\omega_2=50000$ (from the uniform partition), when:
\begin{itemize} 
\item the ratio of $r_k$ is larger than 2;
\item the ratio of the number of iterations is larger than 1.2.
\end{itemize}

\begin{figure}[htbp]
\centering
\includegraphics[angle=-90, width=\linewidth]{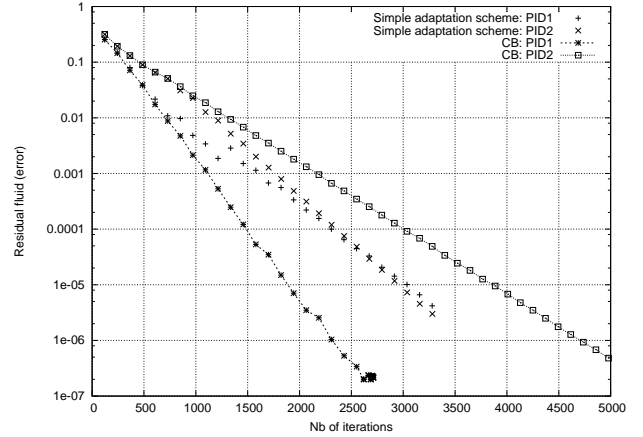}
\caption{Example of adaptation scheme.}
\label{fig:adap}
\end{figure}

Figure \ref{fig:adap-omega} shows the evolution of the $\omega_2$ adaptation.
In this case, the simple adaptation scheme improved very slightly (about 4\%) the CB partition based
cost.

\begin{figure}[htbp]
\centering
\includegraphics[angle=-90, width=\linewidth]{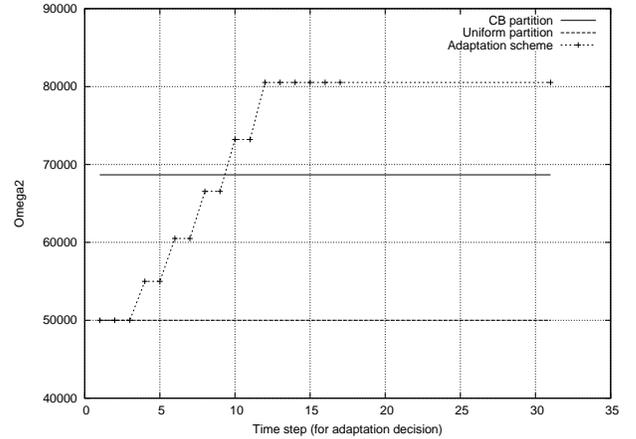}
\caption{Evolution of $\omega_2$.}
\label{fig:adap-omega}
\end{figure}

As it was mentioned earlier, we have to control both the residual fluid (y-axis) and the
computation cost (x-axis). Of course, we have to adapt such ideas to a practical implementation
architecture. And stabilizing such a solution may not be obvious in a general case.\\

The results presented here are yet simulation based evaluation of the computation
efficiency of an asynchronous distributed architecture based on the D-iteration method.
We believe the model used here is realistic enough and shows a promising new
solution to solve linear equations.

\section{Conclusion}\label{sec:conclusion}
In this paper, we presented preliminary simulation results on the evaluation of a distributed
solution associated to the D-iteration method. The main aim of this study was to estimate
the potential of our approach before a real implementation and benchmarking
of such a solution. We showed that our algorithm is very well suited for the parallel
computation of really large linear systems for at least three reasons: firstly, each entity (PID) needs to keep and
update only the local information (distributed memory), the locally kept information
size decreasing linearly with $K$ the number of PIDs; secondly, the efficiency of 
the parallel computation increases with the size $N$ and finally, the computation speed
increases almost linearly with $K$ when $N/K$ is large enough.


\end{psfrags}
\bibliographystyle{abbrv}
\bibliography{sigproc}

\end{document}